\documentclass[11pt,reqno]{amsart}
\usepackage{amsmath,bm}
\usepackage{amsfonts}
\usepackage{amssymb}
\usepackage{amsthm}
\usepackage{amscd}

\usepackage{amsmath}
\usepackage{amssymb}
\usepackage{amsthm}
\usepackage{graphics}
\usepackage{eucal}
\usepackage{mathrsfs}

\usepackage[utf8]{inputenc} % Any characters can be typed directly from the keyboard, eg éçñ
\usepackage{textcomp} % provide lots of new symbols
\usepackage{graphicx}  % Add graphics capabilities
\usepackage{flafter}  % Don't place floats before their definition
\usepackage{bm}  % Define \bm{} to use bold math fonts
\usepackage{tikz}
\usepackage[pdftex,bookmarks,colorlinks,breaklinks]{hyperref}  % PDF hyperlinks, with coloured links
\usepackage{color}
\usepackage{memhfixc}  % remove conflict between the memoir class & hyperref
\usepackage{pdfsync}  % enable tex source and pdf output syncronicity
\usepackage{amssymb}
\usepackage{amsfonts}
\theoremstyle{plain}
\newtheorem*{maintheorem*}{Main Theorem}
\newtheorem*{thm*}{Theorem}
\newtheorem*{thma*}{Theorem A}
\newtheorem*{thmaa*}{Theorem A'}
\newtheorem*{thmb*}{Theorem B}
\newtheorem*{thmo*}{Theorem 1.1}
\newtheorem*{thmc*}{Theorem C}
\newtheorem*{thmd*}{Theorem D}
\newtheorem*{thmf*}{Theorem 4.1}
\newtheorem*{remark*}{Remark}
\newtheorem*{conjecture*}{Conjecture}
\newtheorem*{prop*}{Proposition}
\newtheorem*{lem*}{Basic Lemma}
\newtheorem{thm}{Theorem}[section]

\newtheorem{lem}[thm]{Lemma}

\newtheorem*{proofc*}{Proof of Theorem C}

\begin{document}

%\thanks {\textit{Keywords: Diophantine approximation}}
\author{Youssef Lazar}
%\address{{\tiny Y. Lazar,  Al-Imam Mohammed ibn Saud Islamic University, Riyadh, Saudi Arabia}}
\email{ylazar77@gmail.com}
\date{}

\title{A remark on a conjecture of Erd\H{o}s and Straus}

%%%%%%%%%%%%%%%%%%%%% Abstract %%%%%%%%%%%%%%%%%%%%%%%%%%%%%%%%%%%%%%%%%%%%%%%%%%%%%%%%%%%%%%%%%%%%%%%%%%%
\maketitle

\begin{abstract}  The aim of this note is to show that given a positive integer $n \geq 5$, the positive integral solutions of the diophantine equation $4/n = 1/x + 1/y+1/z$ cannot have solution such that $x$ and $y$ are coprime with $xy < \sqrt{z/2}$. The proof uses the continued fraction expansion of $4/n$.
  \end{abstract}

%%%%%%%%%%%%%%%%%%%%%%%%%%%%%%%%%%%%%%%%%%%%%%%%%%%%%%%%%%%%%%%%%%%%%%   Introduction   %%%%%%%%%%%%%%%%%%%%%%%%%%%%%%%%%%%%%%%%%%%%%%%%%%%%%%%%%%%%%%%%%%%%%%%%%%%

\section{The result}
Given a positive integer $n$ we are interested in the following diophantine equation 
\begin{equation}
\label{egypt1}
\frac{4}{n} = \frac{1}{x}+ \frac{1}{y}+ \frac{1}{z}.
\end{equation}
where $x,y,z$ are three positive integers to be found. The question of finding positive integer solutions $(x,y,z)$ for this equation was raised by  P. Erd\H{o}s and E. Straus (see e.g. \cite{erdos}). This problem has attracted a lot of attention and  despite many efforts it is still widely open. For an account of the main contributions on this subject we refer the reader to \cite{tao}.

\noindent Usual reductions allow us to assume that $n$ is an odd prime number, thus one can assume that $n$ is an odd prime number $ p $. We are looking for integral primitive solutions $ (x,y) \in \mathbb{N} $ of the diophantine problem,
when $z $ is a fixed positive number 
\begin{equation}
\label{egypt2}
 \left|\frac{4}{p} - \frac{x+y}{xy} \right| = \frac{1}{z}.
\end{equation}

\noindent Our main result shows that a solution $(x,y,z)$ of (\ref{egypt2}) cannot have its largest coordinate, say $z$, too far away from the two  other coordinates $x$ and $y$, provided $ \gcd(x,y)=1 $. This gives a certain piece of information regarding the localization of the lattice points which are solution of the problem.  Let us denote by $f(p)$ the number of triples of solutions to (\ref{egypt1}) with $n=p$ or equivalent to (\ref{egypt2}). Recently Elsholtz and Tao gave in \cite{tao} precise bounds for averages of the form $ \sum_{p \leq N} f(p)$, the study of the counting function $f(p)$ reduces to count solutions triples of two different types
\begin{itemize}
\item \textit{Type} I  solutions of (\ref{egypt1}) where  $p$ divides $x$ and $\gcd(p,yz)=1$.
\item \textit{Type} II  solutions of (\ref{egypt1})   where  $p$ divides $y,z$ and $\gcd(p,x)=1$.

\end{itemize}
We denote by $f_{I}(p)$ (resp. $f_{II}(p)$) the number of solutions of (\ref{egypt1}) of type I  (resp. of type II). It was stressed in the same work that for any odd prime number one has the relation 
$$ f(p) =3 f_{I}(p) +3 f_{II}(p).$$
We introduce a new type which we call the type III which are the solutions of (\ref{egypt1}) where  $xy < \sqrt{z/2}$ and $\gcd(x,y)=1$. Analogously we denote by  $f_{III}(p)$ the number of solutions of (\ref{egypt1}) of type III.  Our main result is that for any prime number greater than 3, we have $f_{III}(p)=0$, more precisely
\begin{thm}\label{main}Given an arbitrary prime number $ p \geq 5$, there are no triple of positive integers $ (x,y,z) $  which is solution of (\ref{egypt2}) in the range $\displaystyle xy < \sqrt{z/2}$ and with $ \gcd(x,y)=1 $.
\end{thm}
\noindent \textbf{Proof of the Theorem.} Suppose we fix an arbitrarily large integer $z_{0}>0$ in the range $\displaystyle xy < \sqrt{z/2}$ and let us try to solve the following diophantine equation with $(x,y)\in \mathbb{N}^2$, 

\begin{equation}
\label{egypt3}
 \left|\frac{4}{p} - \frac{x+y}{xy} \right| = \frac{1}{z_{0}}.
\end{equation}

\noindent In particular the equation in (\ref{egypt3}) gives rise to the following inequality 
\begin{equation}\label{ineq1}
\displaystyle \left|\frac{4}{p} - \frac{x+y}{xy} \right| < \frac{1}{2(xy)^2}. 
\end{equation} 
The conclusion of the theorem will follow from the fact that such $x$ and $y$ would never exist.
\noindent Since we assume that $ \gcd(x,y)=1  $ then  we have $ \gcd(x+y,xy)=1$. 
We need the following  classical result of the theory of continued fraction which can be found for instance in \cite{khinchin} (Theorem 19).
\begin{lem}\label{cv} Let $m,n$ be two positive integers and suppose $\gcd(r,s) =1$. If 
$$ \left| \frac{m}{n} - \frac{r}{s} \right| < \frac{1}{2s^2}.$$
Then $\displaystyle \frac{r}{s}$ is one of the convergent of $ \displaystyle \frac{m}{n}$.

\end{lem}

\noindent  By Lemma \ref{cv} we infer from (\ref{ineq1}) that the rational number $\displaystyle \frac{x+y}{xy}$ must be  one of the convergents of $\displaystyle 4/p $. If we write the continued fraction expansion of $\displaystyle 4/p = \left[0; a_1, \ldots, a_l \right]  $, we can say that there exists some $1 \leq k \leq l$ such that 
 $$\displaystyle \frac{x+y}{xy} = c_{k}(\frac{4}{p}) = \left[0; a_1, \ldots, a_k \right] =\frac{p_{k}}{q_{k}}.$$ 
Since $x,y$ play a symmetric role, we can assume that $x \leq y$.  Our fractions are reduced thus we deduce that we might have $ x+y = p_{k}$ and $ xy = q_{k}$. The fact that  such $x$ and $y$ might exist relies on the solvability in $ \mathbb{N} $ of the following quadratic equation 
$$ X^2 - p_k X + q_k = 0.$$
The discriminant $D_{k}=p_k^2 - 4 q_k $ cannot vanish otherwise $ p_k $ and $ q_k$ will fail to be coprime. If $D_{k}=p_k^2 - 4 q_k > 0$, then a couple $(x,y)$ of rational solutions of the equation 

\begin{equation}
\label{egypt4}
\frac{4}{p} = \frac{1}{x}+ \frac{1}{y}+ \frac{1}{z_{0}}
\end{equation}
is given by 
\begin{center}
$\displaystyle x_k = \frac{p_k - \sqrt{D_k}}{2}$ \ and  $\displaystyle y_k= \frac{p_k + \sqrt{D_k}}{2}$. 
\end{center}

\noindent Necessarily  $\displaystyle z_{0}= \frac{p} {4 - p c_{k}} $.
Indeed,  $$\frac{4}{p} -\frac{1}{x_{k}}- \frac{1}{y_{k}} - \frac{1}{z_{0}} = \frac{4}{p} - c_{k} - \frac{4 - p c_{k}}{p} =0.$$  
\noindent  Regarding $x_k$ and $y_k$, both they are in the quadratic field $ \mathbb{Q}[\sqrt{D_{k}}] $, so these are not necessarily integers. In order to obtain integral solutions we are forced to assume that $ D_{k} $ has a square root which is an odd integer. In others words,  $ D_{k} = a^{2}$ where  $ a $ is an odd integer. In this case, we obtain that a triple of solutions  $ (x_{k}, y_{k}, z_{0}) $ which is given by 

\begin{center}
$\displaystyle x_k = \frac{p_k - a}{2}$ \  $\displaystyle y_k= \frac{p_k + a}{2}.$  
\end{center}

\noindent Note that necessarily
\begin{center}
$\displaystyle z_{0}=\left|   \frac{1} {4/p -  c_{k}} \right| = \frac{1}{r_k }$
\end{center}
where $r_k= \left|  4/p -  c_{k} \right| $ is the $k$-th error term in the continued fraction approximation. 
It is well known (see e.g. ) that 
$$ r_k = \frac{1}{q_{k}(x_{k+1} + q_{k-1})} $$

\noindent where $$ x_{k+1} = [a_{k+1}; a_{k+2}, \ldots, a_{l}].$$ 

\noindent Hence, the only possible triples of solutions of (\ref{egypt3}) in the range given above with a fixed $z_0$ must take the following form 
\begin{center}
$\displaystyle x_k = \frac{p_k - a}{2}$ \  $\displaystyle y_k= \frac{p_k + a}{2}$  and $ z_{0} = q_{k}(x_{k+1}q_{k} + q_{k-1})$
\end{center}
where $\displaystyle \frac{p_{k}}{q_{k}}$ is one of the convergents of the continued fraction expansion $ 4/ p = [a_{1}; a_{2}, \ldots, a_{l}] $ and provided $ p_{k}^{2}-4q_{k} = a^2$ with $a$ being an odd integer. We will show that the latter condition can never be fullfilled. \\
To proceed we take advantage from the fact  that the convergents of $4/p $ can only assume specific values which are given in the following lemma, 
\begin{lem}\label{cv2} For any prime number $ p \geq 5 $ set \begin{center}
$$a_{1}= \left\lbrace \begin{array}{cc}
\displaystyle   \frac{p-1}{4}  &  \ \mathrm{if} \ \ \ p \equiv 1 \pmod{4} \\
 \displaystyle   \frac{p -3}{4}  &  \ \mathrm{if} \ \ \ p \equiv 3 \pmod{4}.
\end{array}\right.$$

\end{center} We have two cases, 
\begin{itemize}
\item[$(a)$] if $ p \equiv 1 \pmod{4} $, then  $ \displaystyle \frac{4}{p} = [0 ; a_{1},4] $ and the convergents are $\displaystyle   \{0,  \frac{4}{p-1},   \frac{4}{p}    \}. $
\item[$(b)$] If $ p \equiv 3 \pmod{4} $, then $ \displaystyle \frac{4}{p} = [0 ; a_{1},1,3] $ and the convergents are  $\displaystyle   \{0,  \frac{4}{p-3},   \frac{4}{p+1},   \frac{4}{p}    \}. $

\end{itemize}

\end{lem}
\noindent \textit{Proof.} The continued fraction of a rational number is entirely determine by the eucliden algorithm between $4$ and $p$. 

 $(a)$ Suppose $ p \equiv 1 \pmod{4} $, we perform to the division algorithm
 $$ 4 = p (0) + 4$$
  $$ p = 4 (a_{1}) + 1$$
  $$ 4 = 1 (4) +~ 0.$$
Here  $\displaystyle a_{1}= \lfloor \frac{p}{4} \rfloor = \frac{p-1}{4}$, thus $\displaystyle \frac{4}{p} = [0 ; a_{1},4]  $. The successive convergents are given by 
$ c_{0}=0 $, $ \displaystyle  c_{1}= \frac{1}{a_{1}}= \frac{4}{p-1} $ and $\displaystyle  c_{2} =  \frac{1}{a_{1} + \frac{1}{4}}= \frac{4}{p}$.

 $(b)$ Suppose $ p \equiv 3 \pmod{4} $, the division algorithm again shows that
 $$ 4 = p (0) + 4$$
  $$ p = 4 (a_{1}) + 3$$
  $$ 4 = 3(1) +~ 1$$
  $$ 3 = 1(3) +~ 0.$$
  
\noindent   Therefore $\displaystyle \frac{4}{p} = [0 ; a_{1},1,3]  $.  The corresponding convergents are $c_{0}=0$, $c_{1}=\frac{4}{p-3}$,\\ $\displaystyle c_{2}=\frac{1}{a_{1}+1}= \frac{1}{\frac{p-3}{4}+1}= \frac{4}{p+1}$ and $\displaystyle c_{3}=\frac{1}{a_{1}+ \frac{1}{1+ \frac{1}{3}}}= \frac{1}{a_{1}+ \frac{3}{4}}= \frac{4}{p}$.
This proves the Lemma.\\

We are ready to conclude. In both case, Lemma \ref{cv2} shows that all the non-trivial convergents of $4/p$ (i.e. other that $0$ and $4/p$) are egyptian fractions, in particular $p_{1}=1$ (case $(a)$) and $p_{1}=p_{2}=1$ (case (b)).  It follows that the $X^2 - p_{k} X + q_{k} =0$ is not solvable in $\mathbb{N}$ since $D_{k}=1-4q_{k} < 0$ in all the cases. Hence in the range $\displaystyle xy < \sqrt{z/2}$  there are no solution of (\ref{egypt2}) with $ \gcd(x,y)=1$. This finishes the proof of Theorem \ref{main}.
\section{Concluding remarks}
\noindent  Let $p$ be a prime number, and let us introduce the subset 
$$E_{p} : = \left\lbrace (x,y,z) \in \mathbb{N}^{3} ~~\vert~~    4/p = 1/x+1/y + 1/z  \right\rbrace.$$ Then
$p$ is solution of (\ref{egypt2}) if and only if $E_{p} \neq \emptyset$. So the validity of the Erdos-Straus conjecture amounts to prove that $E_{p} \neq \emptyset$ for every prime $p$. The Theorem \ref{main} tells us that for any prime $ p \geq 5 $ the set

$$R_{p} : = E_{p} \cap \left\lbrace (x,y,z) \in \mathbb{N}^{3} ~~\vert~~  ~~\gcd(x,y)=1, \displaystyle xy < \sqrt{z/2} \right\rbrace$$
is empty, in other words $f_{III}(p)=0$. Now suppose we are looking for solutions of (\ref{egypt1}) with $0 <x,y,z \leq N$, so we have a total number of possibilities equal to $(N-1)^3$ from which we have to remove the elements of $ \left\lbrace (x,y,z) \in [1,N]^3~~\vert~, ~~\gcd(x,y)=1, \displaystyle xy < \sqrt{z/2} \right\rbrace$. Let us denote by $A_{N}$ this subset and put $a_N = |A_N|$. An estimate of $ a_{N} $ can be obtained. Indeed by slicing we get,
$$  a_{N}= \sum_{1\leq z\leq N}|\left\lbrace 1 \leq x,y \leq  \lfloor \sqrt{z/2} \rfloor~~\vert~, ~~\gcd(x,y)=1, \displaystyle xy < \sqrt{z/2} \right\rbrace \vert.$$
The inner sum counts the number of primitive lattice points under the hyperbola of equation $\displaystyle y= \sqrt{z/2} ~ x^{-1}$. Thus we can write 
$$  a_{N}= \sum_{1\leq z\leq N} \sum_{x \leq  \lfloor \sqrt{z/2} \rfloor~}\sum_{xy \leq  \lfloor \sqrt{z/2} \rfloor~}\varphi(y). $$
We need the following estimate which can be found in \cite{tene} (Thm 3.4) as $X \rightarrow \infty$
$$\sum_{n \leq  X~}\varphi(n) = \frac{6}{\pi^2} X^2 + O(X\ln X).$$

\noindent This asymptotics leads us to 
$$  a_{N}=\frac{3}{\pi^2}  \sum_{1\leq z\leq N} \sum_{x \leq  \lfloor \sqrt{z/2} \rfloor~}  \frac{z}{x^2} + O( \frac{\sqrt{z}}{2x}( \ln \frac{z^{1/2}}{x})). $$

\noindent  \noindent Using crude bounds we get
$$   \frac{3}{\pi^{2}} N(1 + o(1))   \ll    a_{N} \ll \frac{1}{2\sqrt{2}\pi^{2}} N^{5/2} + O(N^{3/2}\ln N).$$
As a conclusion,  the previous discussion shows that the number of lattices points to be discarded from $ E_{p} \cap [1,N]^3$ is at an order of magnitude of at least $ N $ and at most $ N^{5/2} $ lattices points inside the cube $ [1,N]^3 $. \\

\textit{Further questions.} It would be interesting to find an analog of the main result for $ 5/p $ instead of  $ 4/p $, which is also a conjecture due to Sierpinski. More generally, one can also try to ckeck what is happening for fractions of the form $a/p$ where $ a $ is a integer which does not divide $ p $.  The continued fraction expansion will depend on the class of $ a \mod p $ and in this case Lemma \ref{cv2} might be less trivial than the $ 4/p $-case.

\end{document}